\documentclass[12pt]{article}
\usepackage{amsmath}
\usepackage{epsfig}
\font\Bbb=msbm10 at 11pt  
\begin{document}
\begin{center}
{\bf The Elliptical Envelope}
\vspace{11pt}

Stephen Howard

\vspace{11pt}
Copyright of the Mathematical Association of America 2007. All rights reserved.
\end{center}

\paragraph{Introduction.}

An \emph{isocline} of a family of curves is the locus of all points in the family having the same specified slope \cite[p. 15]{NagleSaff}. 

\emph{Problem:} Under what circumstances can an isocline of a family of curves also be a member of the family? 

While this remains an open problem in general, one interesting and highly symmetric special case involves a one-parameter family of circles in ${\hbox{\Bbb R}}^2$ that I will refer to as $F_C$, which does have this property. It can be uniquely generated by specifying that its zero slope isocline (or \emph{null isocline}) is a unit circle. If each circle $C_{\alpha}$ in the family has a size and position that is constrained by the requirement that the null isocline $C_{0}$ intersects $C_{\alpha}$ at the top and bottom of $C_{\alpha}$, (the points of maximum and minimum $y$-values, which are the points of zero slope), then the family of circles is $F_C$.	Each circle in $F_C$ is centered along the $x$-axis at some displacement $\alpha$. The radius of the circle centered at $(\alpha,0)$ must be equal to the height of the null isocline at $x = \alpha$, so  \(r(\alpha) = \sqrt{1-\alpha^2}. \) Then $F_C$ can be parameterized by
\begin{equation}
(x- \alpha)^2 +y^2 = 1-\alpha^2,    \label{equ_FC} 
\end{equation}
with $\alpha$ varying between $-1$ and 1.
\begin{figure}[!ht]
\begin{center}
\epsfig{figure=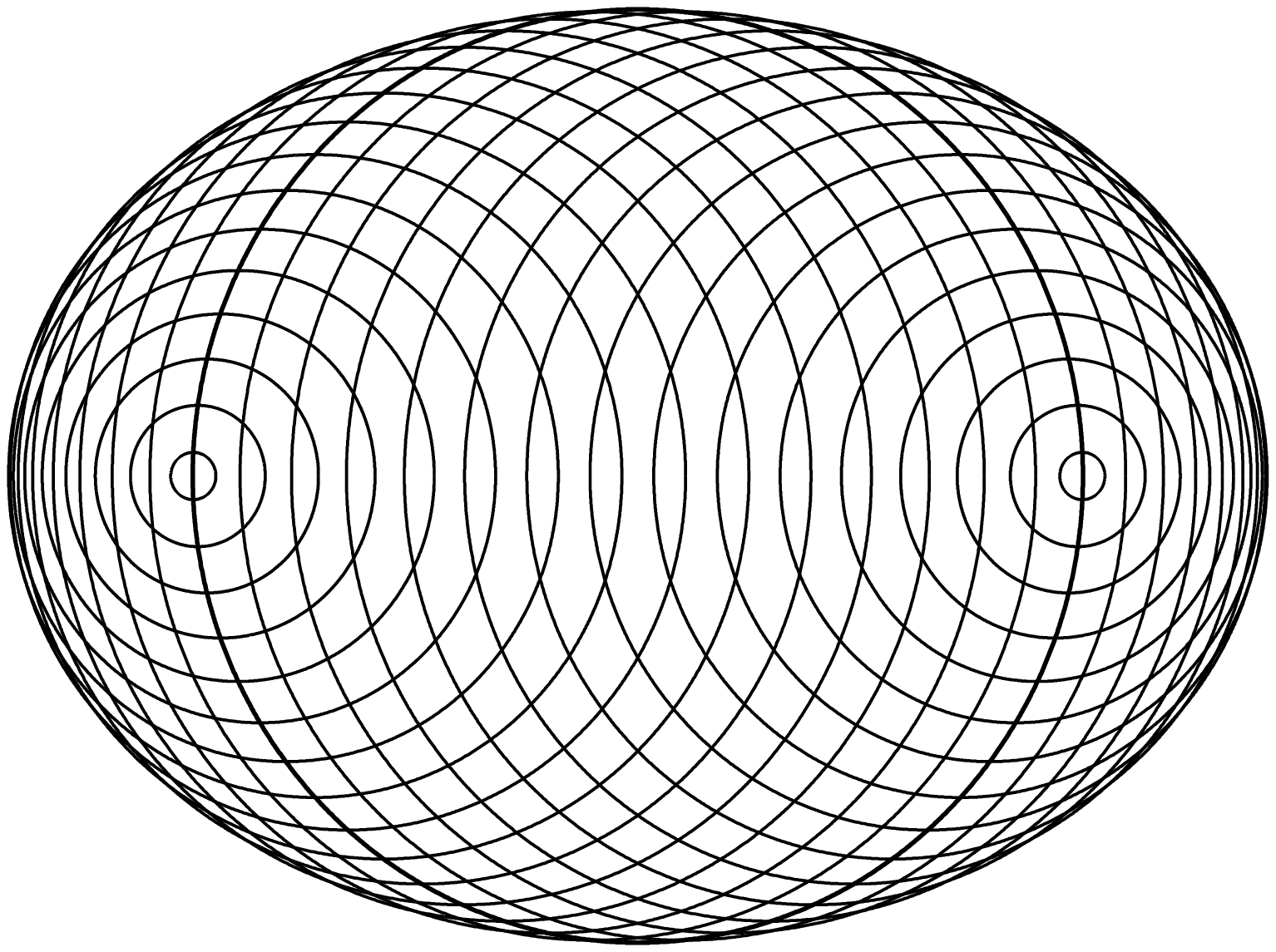, height=120pt, width=169.2pt}
\caption{Family of circles $F_C$}\label{fig:FC}
\end{center}
\end{figure}

The most basic observation that can be made about this family of curves is that it is bounded. What is the shape of this boundary or envelope? 

It appears to be an ellipse. Not all circles in the family contribute points to the boundary curve, and consideration of the outermost circle to do so yields the fact that it has the same major and minor axes as the ellipse, 
\begin{equation}
\frac{1}{2} {x^2} + y^2 = 1.  \label{equ_EB}
\end{equation}  

However, without further analysis we can not rule out the possibility that the boundary could take any ellipse-like path between the points \((0, \pm 1) \) and \( (\pm \sqrt 2 ,0) \), without being truly elliptical. 

To resolve this we need to precisely define the boundary and derive the equation that it obeys. This can be done in a variety of ways. The classical method is based on the fact that the envelope is locally tangent to the family of curves, and by applying the \emph{fundamental theorem of envelopes} \cite{Biswas:mmt} we can derive an equation for the boundary by finding the roots of the partial derivative with respect to $\alpha$ of the \emph{level set} formula for the family $F_C$. This proceeds rapidly,
\begin{eqnarray}
	f(x,y,\alpha) &=& (x- \alpha)^2 +y^2 +\alpha^2 - 1 = 0, \nonumber\\ 
	{\frac{\partial{f(x,y,\alpha)}}{\partial{\alpha}}} &=& -2(x- \alpha) + 2\alpha = -2x + 4\alpha = 0.\nonumber
\end{eqnarray}
Hence the boundary point obeys $x_b = 2 \alpha$, and when this is substituted back into the parametric formula, we find the elliptical relation for the shape of the envelope (\ref{equ_EB}).

Some important and non-obvious details of the validity of this classical theory of envelopes have been glossed over here, but there are outside the main purpose of this article \cite{AndersKock}, \cite{Kunihiro:ptp}.  Rather, this family of circles is really more interesting as a testbed for devising ways to make the subtle, infinitesimal nature of the envelope become more comprehensible and concrete. In this article we will consider three nonstandard methods, each with a distinct definition for the envelope of $F_C$, to illustrate the relation between how elegantly you define the problem and how easily it is solved. 
The methods are presented approximately in order of decreasing obviousness, and increasing elegance (in my opinion).   

\paragraph{Boundary Definition 1:}  In the first method, we use polar coordinates $(r,\theta)$ and define the boundary to be the set of all points $\left(R(\theta),\theta \right)$ that are maximally far away from the origin 
\[ 			B=\{(R(\theta),\theta)\in F_C \mid\forall (r(\theta),\theta) \in F_C,\ r(\theta) \le R(\theta)\}.\]

Each point in the boundary set $B$ lies on an element of $F_C$ and has a unique $\theta$ coordinate. The boundary point in the $\theta$  direction is the point within $F_C$ of maximum radial coordinate in that direction. So to find the boundary point in the $\theta$ direction, we can find the intersection of a ray emanating from the origin with each circle that it crosses. Then we find the intersection point that is furthest away from the origin. This method is at least conceptually straightforward and was the first to be considered when working this problem, by myself as well as by other undergraduate students when they independently tried to solve the envelope problem \cite{howard97:hrumc}. 
\begin{figure}[!ht]
\begin{center}
\epsfig{figure=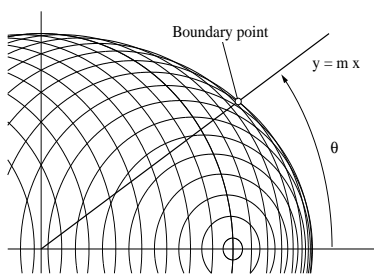, height=120pt, width=169.2pt}
\caption{Boundary point defined by Method 1}
\end{center}
\end{figure}

So we consider a ray oriented  at an angle $\theta$ with respect to the x-axis.
\begin{equation}
     y = mx,   \quad \quad m = tan(\theta)     
\end{equation}
From here on it is convenient to work through the derivation in $(x,y)$ Cartesian coordinates. Substitution of (3) into (\ref{equ_FC}) yields the intersection of this ray with circle $C_\alpha$ 
\[(x-\alpha)^2 +(mx)^2 = 1-\alpha^2, \]
or
\[(m^2+1)x^2 -2\alpha x +2\alpha^2 -1 = 0. \]
Solved for $x$ this becomes
\begin{equation}   
x(\alpha,m) = {\frac{\alpha \pm \sqrt{\alpha^2-(m^2+1)(2\alpha^2 -1)}}{m^2+1} }, 
\end{equation}
and gives the intersection point \((x(\alpha,m),m x(\alpha,m))\) of the $\alpha$-circle with the $m$-ray. The distance from the origin to the intersection point is 
\[\sqrt{(x^2 +(mx)^2)} \quad {\mbox{or simply}} \quad |x|\sqrt{(1+m^2)},\]
thus it is maximized for a given $m$ when $|x|$ is maximized. 

To find the boundary point at $m$ we simply find what $\alpha$ maximizes $|x|$, which occurs when $dx/d\alpha = 0$. The result of the differentiation is:
\begin{equation}
{\frac{dx}{d\alpha}} = {\frac{1}{m^2+1}}\Bigg {(}1 \pm {\frac{(2m^2+1)\alpha}{\sqrt{\alpha^2-(m^2+1)(2\alpha^2 -1)}}}\Bigg {)},      
\end{equation}
and this equals zero when        
\[2(2m^2 +1)^2\alpha^2 =  m^2 +1.\]
Solved for $\alpha$ this becomes
\begin{equation}
\alpha_{max} = {\frac{1}{2\sqrt{(m^2 +1/2)}}}, \label{equ:amax} 
\end{equation}
the circle at $\alpha_{max}$ contains the boundary point at the angle $\theta = arctan(m)$. To locate this point we will plug (6) into (4) and simplify,
\begin{equation}
x_b(m)  = x(\alpha_{max}, m)= \frac{1}{\sqrt{(m^2+1/2)}}. \label{equ:xb} 
\end{equation}
It is now possible to show that the boundary is elliptical. We can rephrase the last result as
\[m^2 {x_b}^2 + \frac{1}{2} {x_b}^2  = 1, \]
then substitute $y_b = mx_b$ and we have derived the equation of the boundary: 
\[ {y_b}^2 + \frac{1}{2} {x_b}^2  = 1. \]   
And it is indeed exactly elliptical, in agreement with our original guess.\\    

It is interesting to compare (\ref{equ:amax}) and (\ref{equ:xb}). We see again the same initial result we found by the standard method; that a circle centered at $\alpha$ always has a boundary point at $x_b = 2\alpha$, if that circle has a boundary point at all. A circle will only contribute a point to the boundary if it is centered at $\alpha \le \sqrt{2}/2$.

\paragraph{Method 2.} The second, faster method involves directly taking a limit to calculate the boundary points instead of using a differential maximization \mbox{\cite[pp. 169--179]{Courant:difint}}. Now, the boundary point contributed by a circle $C_\alpha$ centered at $(\alpha,0)$ is defined by finding the intersection point of $C_\alpha$ with a second circle $C_{(\alpha+\delta)}$ centered at $(\alpha+\delta,0)$ and then taking the limit of the intersection as $\delta \rightarrow 0$.

\paragraph{Boundary Definition 2:}  Let $B_\alpha$ be the boundary points contributed to the boundary by circle $C_\alpha$. These boundary points are the limit of the points of intersection with neighboring circles,
\[B_\alpha \equiv {\lim_{\delta \rightarrow 0}} (C_\alpha \cap C_{(\alpha+\delta)}) \]  
when this limit exists (it does not when $\alpha > \sqrt{2}/2$). The boundary can be arbitrarily well approximated using a finite 
number of circles in the family and joining the exterior arcs of the circles between the points of intersection with adjacent circles. In fact figure \ref{fig:FC} depicting $F_C$ uses only 31 circles and does a nice job of approximating the envelope. In the limit as $\delta \rightarrow 0$ (and number of circles $\rightarrow \infty$), the exterior arcs will approach zero length and adjacent intersection points coincide with each other, and lie exactly on the boundary.
\begin{figure}[hb]
\begin{center}
\epsfig{figure=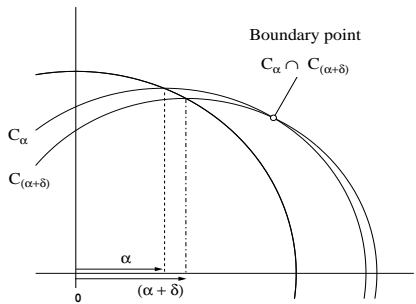, height=120pt, width=169.2pt}
\caption{Boundary point defined by Method 2}\label{fig:BM2}
\end{center}
\end{figure}
To proceed using this definition, the boundary point contributed by $C_{\alpha}$ can be solved for, by finding the intersection points of $C_{\alpha}$ and $C_{(\alpha+\delta)}$. We solve the system
\begin{eqnarray}
	C_{\alpha}{\!\!}&:&    \quad		(x -\alpha)^2 + y^2 = 1 - {\alpha^2},  \nonumber  \\                     
	C_{(\alpha+\delta)}{\!\!}&:& \quad	(x - \alpha -\delta)^2 +y^2 = 1 - (\alpha+\delta)^2.  \nonumber
\end{eqnarray}
Eliminating y gives
\[(x-\alpha)^2 +\alpha^2 - (x-\alpha-\delta)^2 - (\alpha+\delta)^2 = 0.\]
Expanding this,
\[(x^2 -2 \alpha x+ \alpha^2) + \alpha^2 -(x^2-2 \alpha x - 2 \delta x +\alpha^2 + 2 \alpha \delta +\delta ^2) -(\alpha^2+2 \alpha \delta +\delta ^2) = 0, \] 
then collecting like terms
\[2 \delta x - 4\alpha \delta - 2\delta ^2  =  0,\]
divide by $2\delta$ and we have
\[x = 2\alpha + \delta.  \]
This is the $x$-coordinate of the intersection point of  $C_{\alpha}$ and $C_{(\alpha+\delta)}$. Taking the limit as $\delta \rightarrow 0$  gives the boundary point,
\[ x_b  =  \lim_{\delta \rightarrow 0} x  =  2\alpha, \]
which agrees with the result we observed before. So we have found that $x_b = 2\alpha$ or $\alpha = \frac{1}{2} x_b$ plugging this into (\ref{equ_FC}) to find the $y$-coordinate of the boundary point in terms of the $x$-coordinate, we have
\[(x_b - \frac{1}{2}x_b)^2 + y_b^2  = 1-(\frac{1}{2}x_b)^2,\]
\[\frac{1}{4}x_b^2 +y_b^2 = 1 - \frac{1}{4}x_b^2,\]
\[\frac{1}{2}x_b^2 + y_b^2 = 1,\]
which is again the equation of the elliptical envelope.\\

\begin{figure}[!ht]
\begin{center}
\epsfig{figure=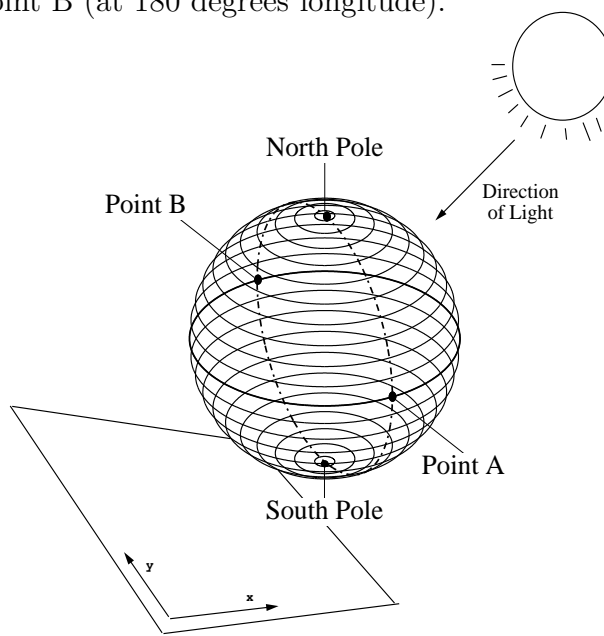, height=230pt, width=230pt}
\caption{Method 3, Yoga ball and surface} \label{fig:yogaball}
\end{center}
\end{figure}

\paragraph{Method 3: Proof by Projection.} There exists a higher-dimensional object that leads to $F_C$ by way of projection. Consider a transparent sphere of unit radius. For the purposes of fixing a definite orientation in space we will designate North and South poles of the sphere, and then to begin the proof we will draw some circles on the sphere's surface. In particular we will consider the complete family of latitudinal circles, including the equator, which are conveniently visualized via a finite number these circles drawn evenly over both hemispheres. In addition, one circle of longitude is drawn passing through N and S poles. A sphere labeled in such a way looks a lot like a transparent yoga exercise ball, which was the seed of inspiration for this method. Refer to figure \ref{fig:yogaball}.  

Pick one of the intersection points of the equator and the longitudinal circle, and call that point A (defined to be zero degrees longitude, zero degrees latitude). Then the other intersection point on the other side of the sphere can be called point B (at 180 degrees longitude).

Now consider this sphere positioned in space above a plane, so that the center of the sphere is directly above the origin of a Cartesian $(x,y)$ coordinate system on the plane. Orient the sphere to make the diameter going through A B (the A-B axis) become parallel to the $y$-axis of the plane. We will allow the sphere to be rotated about the A-B axis, but the orientation of the axis itself is fixed.

Next, shine light perpendicular to the plane, passing through the sphere so that shadows of the circles will be cast onto the plane. Since the circles may be oriented at some angle $\theta$ relative to the plane, the shadows projected onto the plane will have an elliptical shape. The major radius of each ellipse will be the same as the radius of the corresponding circle on the sphere, while the minor radius of the ellipse will be shortened by a factor of $\cos \theta$. If we then set $\theta = \pi /4$ , the shrinking factor in the $x$ direction will be $\frac{1}{\sqrt{2}}$, and the shadow of the equator will exactly overlay the shadow of the longitudinal circle. More importantly, the set of  elliptical projections on the surface will display a strong resemblance to the set of circles $F_C$ that we were originally considering. Let us call the family of elliptical projections $F_E$.

\begin{figure}[!ht]
\begin{center}
\epsfig{figure=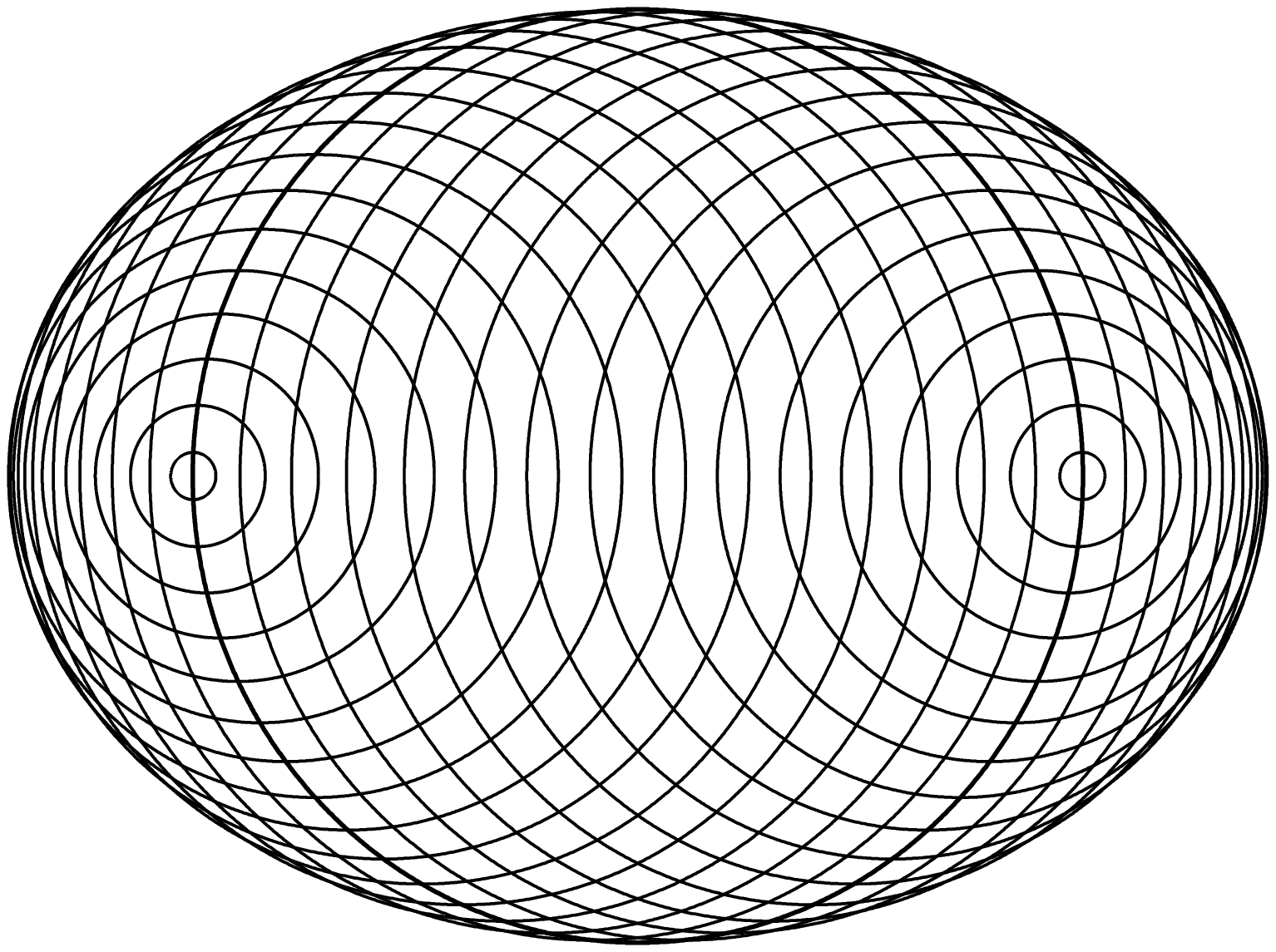, height=150pt, width=150pt}
\caption{$F_E$ as it appears projected onto the plane}
\end{center}
\end{figure}

The most important similarity between $F_C$ and $F_E$ is that in both families the largest element is the null isocline of the family, and in that way, the largest element generates all the other elements. In $F_E$ this property is guaranteed because the largest element is the projection of the longitudinal circle, which intersects all the other circles on the sphere at 0 degrees longitude, which corresponds to the point of maximum $y$ value of each circle's projection, and also at 180 degrees longitude on the sphere, which gets mapped to the point of minimum $y$ value for each projection. And of course, at the points of maximum and minimum $y$ value we have $dy/dx = 0$, so the projection of the longitudinal circle is the null isocline of $F_E$.  We can then formulate $F_E$ in terms of a parameter $\alpha$. The ellipse in $F_E$ that is centered at $x = \alpha$ is described by: 
\begin{equation} 
2 (x- \alpha)^2 + y^2 = 1- 2 \alpha ^2. \label{equ_FE}
\end{equation}
This construction of $F_E$ by projecting the sphere onto a surface makes it obvious that its boundary curve is simply a circle \(x^2 + y^2 = 1 \). Now the original problem, ``what is the boundary of $F_C$?'', becomes truly transparent. Simply stretch $F_E$ linearly in the $x$ direction by a factor of $\sqrt{2}$ and it becomes congruent to $F_C$. Then the boundary of of $F_C$ is simply the stretched version of the boundary of $F_E$. This linear transformation is $x \rightarrow x' = x \sqrt{2}, \ \alpha \rightarrow \alpha' = \alpha \sqrt{2}$,  so we can make the substitution $x =x' / \sqrt{2},\ \alpha =\alpha' / \sqrt{2} $ into (\ref{equ_FE}) and the result is the relation for the circles in $F_C$
\begin{equation}
(x'- \alpha')^2 +y^2 = 1-(\alpha')^2.    \label{equ_FCprime} 
\end{equation}
The congruency between $F_C$ and the stretched version of $F_E$ is now proved. And so this stretching takes the circular envelope of $F_E$ and turns it into the elliptical envelope of $F_C$, $\frac{1}{2} {(x')^2} + y^2 = 1 $.\\

\paragraph{Discussion.} It is easy to see the difference in the amount of work needed to execute these strategies. Method 1 takes more than three pages to work out completely by hand, requires a differentiation and carrying the slope $m$ through the algebra, while Method 2 takes one page to work out completely by hand and requires no differentiation of radicals. We characterize Method 1 as being essentially extrinsic in that it relies on adding an external construct (the ray) with which to determine the boundary points, while Method 2 and the classical theory of envelopes are entirely intrinsic, using only the relation between the circles to determine the shape of the envelope. Therefore, this is a good example of the power of intrinsic techniques. 

Methods 2 is in a sense more elementary than the standard method, but essentially they follow the same algebraic strategy. Method 2, without acknowledging it as so, makes a direct calculation of the limit definition of the derivative, while the standard approach makes use of the derivative power rule to directly evaluate it, resulting in only 2 lines of derivation supposing that the main result of the theory of envelopes is already in hand. The thought-provoking concept here is the equivalence between the limit of intersection of neighboring curves, the tangent curve, and the envelope.

Alternatively, Method 3, is a good example of the power of working in a higher-dimensional space \cite{DUGUNDJI:IEEETIT}. It is truly satisfying to see that the family $F_C$ in the plane of ${\hbox{\Bbb R}}^2$ is really just the projection of another simpler object in ${\hbox{\Bbb R}}^3$.


\begin{thebibliography}{1}

\bibitem{NagleSaff}
Nagle and Saff.
\newblock {\em Fundamentals of Differential Equations}.
\newblock Addison-Wesley.
\newblock 3rd Edition.

\bibitem{Biswas:mmt}
Arunava Biswas, Michael Stevens, and Gary~L. Kinzel.
\newblock A comparison of approximate methods for the analytical determination
  of profiles for disk cams with roller followers.
\newblock {\em Mechanism and Machine Theory}, 39(6):645--656, June 2004.
\newblock A good example of how the theory of envelopes can be applied.

\bibitem{AndersKock}
Anders Kock.
\newblock {\em Envelopes - notion and definiteness}.
\newblock University of Aarhus, December 2003.
\newblock Preprint, ISSN: 1397--4076.

\bibitem{Kunihiro:ptp}
Teiji Kunihiro.
\newblock A geometrical formulation of the renormalization group method for
  global analysis.
\newblock {\em Prog. Theor. Phys}, 94:503--514, 1995.
\newblock Erratum-ibid. 95 (1996) 835.

\bibitem{howard97:hrumc}
S.~Howard.
\newblock Two solutions to the elliptical envelope problem.
\newblock In {\em Hudson River Undergraduate Mathematics Conference}. Williams
  College, 1997.

\bibitem{Courant:difint}
R.~Courant.
\newblock {\em Differential and integral calculus, Vol II}.
\newblock Interscience Publishers, 1947.
\newblock translated by E.J. McShane.

\bibitem{DUGUNDJI:IEEETIT}
J.~Dungundji.
\newblock Envelopes and pre-envelopes of real waveforms.
\newblock {\em IEEE Transactions on Information Theory}, 4(1):53--57, March
  1958.
\newblock An example of deriving an envelope from a higher dimensional object,
  in this case a complex valued function yields the envelope of a real valued
  family of curves by taking its absolute value.

\end{thebibliography}









\end{document}